\documentclass[psamsfonts,reqno]{amsart}
\usepackage{amssymb, eucal}

\theoremstyle{plain}

\newtheorem{lemma}{Lemma}[section]
\newtheorem{theorem}[lemma]{Theorem}
\newtheorem{proposition}[lemma]{Proposition}
\newtheorem{corollary}[lemma]{Corollary}
\newtheorem*{thm}{Theorem}
\newtheorem*{CLPS}{CLP (semilattice formulation)}
\newtheorem*{KFST}{Kuratowski's Free Set Theorem}
\newtheorem*{sclaim}{Claim}

\theoremstyle{definition}
\newtheorem{definition}[lemma]{Definition}

\newtheorem{problem}{Problem}

\theoremstyle{remark}
\newtheorem{remark}[lemma]{Remark}

\newcommand{\qedsc}{{\qed}~{\rm Claim.}}

\newenvironment{scproof}
{\begin{proof}[Proof of Claim.]}
{\qedsc\renewcommand{\qed}{}\end{proof}}

\numberwithin{equation}{section}

\newcommand{\pup}[1]{\textup{(}{#1}\textup{)}}

\newcommand{\contr}{a contradiction}

\newcommand{\set}[1]{\{#1\}}
\newcommand{\setm}[2]{\set{#1\mid#2}}
\newcommand{\Set}[1]{\left\{#1\right\}}

\newcommand{\seq}[1]{({#1})}
\newcommand{\seqm}[2]{\seq{{#1}\mid{#2}}}
\newcommand{\Seq}[1]{\left({#1}\right)}

\newcommand{\famm}[2]{(#1\mid#2)}
\newcommand{\Famm}[2]{\left(#1\mid#2\right)}
\newcommand{\bt}{\mathbin{\bowtie}}
\DeclareMathOperator{\im}{im}

\newcommand{\Pow}{\mathfrak{P}}

\DeclareMathOperator{\supp}{supp}

\DeclareMathOperator{\rk}{rk}

\newcommand{\cC}{\mathcal{C}}

\newcommand{\cR}{\mathcal{R}}

\newcommand{\cL}{\mathcal{L}}
\newcommand{\cD}{\mathcal{D}}
\newcommand{\cG}{\mathcal{G}}
\newcommand{\es}{\varnothing}

\newcommand{\onto}{\twoheadrightarrow}
\newcommand{\zero}{\bar{0}}
\newcommand{\one}{\bar{1}}

\newcommand{\eps}{\varepsilon}
\DeclareMathOperator{\Sem}{Sem}
\newcommand{\CLR}{\mathrm{CLR}}

\newcommand{\go}{\omega}
\newcommand{\gos}{\go\setminus\set{0}}

\newcommand{\ba}{\boldsymbol{a}}
\newcommand{\bb}{\boldsymbol{b}}
\newcommand{\bc}{\boldsymbol{c}}
\newcommand{\be}{\boldsymbol{e}}
\newcommand{\bx}{\boldsymbol{x}}
\newcommand{\by}{\boldsymbol{y}}
\newcommand{\bz}{\boldsymbol{z}}
\newcommand{\bu}{\boldsymbol{u}}
\newcommand{\bv}{\boldsymbol{v}}
\newcommand{\bw}{\boldsymbol{w}}

\newcommand{\dxi}{\dot{\xi}}

\DeclareMathOperator{\Con}{Con}
\DeclareMathOperator{\Conc}{Con_c}
\newcommand{\Concs}{\Con_{\mathrm{c}}}

\newcommand{\ol}[1]{\overline{{#1}}}

\newcommand{\jz}{$(\vee,0)$}
\newcommand{\ju}{$(\vee,1)$}
\newcommand{\jzu}{$(\vee,0,1)$}
\newcommand{\jzs}{\jz-semi\-lat\-tice}
\newcommand{\jus}{\ju-semi\-lat\-tice}
\newcommand{\jzus}{\jzu-semi\-lat\-tice}
\newcommand{\jzh}{\jz-ho\-mo\-mor\-phism}

\newcommand{\jzuh}{\jzu-ho\-mo\-mor\-phism}
\newcommand{\jze}{\jz-em\-bed\-ding}

\newcommand{\js}{join-sem\-i\-lat\-tice}
\newcommand{\jh}{join-ho\-mo\-mor\-phism}

\begin{document}

\title[The Congruence Lattice Problem]%
{A solution to Dilworth's Congruence Lattice Problem}

\author[F.~Wehrung]{Friedrich Wehrung}
\address{LMNO, CNRS UMR 6139\\
D\'epartement de Math\'ematiques, BP 5186\\
Universit\'e de Caen, Campus 2\\
14032 Caen cedex\\
France}
\email{wehrung@math.unicaen.fr}
\urladdr{http://www.math.unicaen.fr/\~{}wehrung}

\keywords{Lattice; algebraic; distributive; congruence; congruence-compatible; semilattice; weakly distributive; free set}
\subjclass[2000]{06B15, 06B10, 06A12. Secondary 08A30, 08B10}

\dedicatory{Dedicated to George Gr\"atzer and Tam\'as Schmidt}

\date{\today}

\begin{abstract}
We construct an algebraic distributive lattice $D$ that is not isomorphic
to the congruence lattice of any lattice. This solves a long-standing open problem, traditionally attributed to R.\,P. Dilworth, from the forties. The lattice~$D$ has compact top element and $\aleph_{\go+1}$ compact elements. Our results extend to any algebra possessing a congruence-compatible structure of a \js\ with a largest element.
\end{abstract}

\maketitle

\section{Introduction}\label{S:Intro}

For an \emph{algebra} $L$ (i.e., a nonempty set with a collection of operations from finite powers of $L$ to $L$), a \emph{congruence} of~$L$ is an equivalence relation on~$L$ compatible with all operations of~$L$. A map $f\colon L^n\to L$ (for some positive integer~$n$) is \emph{congruence-compatible}, if every congruence of~$L$ is a congruence for~$f$. (This occurs, for example, in case~$f$ is a \emph{polynomial} of~$L$, that is, a composition of basic operations of~$L$, allowing elements of~$L$ as parameters.) For elements $x,y\in L$, we denote by $\Theta_L(x,y)$ the least congruence that identifies~$x$ with~$y$, and we call the finite joins of such congruences \emph{finitely generated}. We denote by $\Con L$ (resp., $\Conc L$) the lattice (resp., \jzs) of all congruences (resp., finitely generated congruences) of $L$ under inclusion. A homomorphism of \js s $\mu\colon S\to T$ is \emph{weakly distributive} at an element $\bx$ of $S$, if for all $\by_0,\by_1\in T$ such that
$\mu(\bx)\leq\by_0\vee\by_1$, there are $\bx_0,\bx_1\in S$ such that
$\bx\leq\bx_0\vee\bx_1$ and $\mu(\bx_i)\leq\by_i$, for all $i<2$. We say that~$\mu$ is \emph{weakly distributive}, if it is weakly distributive at every element of~$S$. (In case both~$S$ and~$T$ are distributive, this is equivalent to the definition presented in~\cite{WeURP}. Moreover, it extends the original definition given by Schmidt \cite{Schm81,Schm82}.)

In the present paper we prove the following result (cf. Theorem~\ref{T:NonRepr}).

\begin{thm}
There exists a distributive \jzus\ $S$ such that for any algebra $L$ with a congruence-compatible structure of a \jus, there exists no weakly distributive \jzh\ $\mu\colon\Conc L\to S$ with~$1$ in its range. Furthermore, $S$ has $\aleph_{\go+1}$ elements.
\end{thm}

As every isomorphism is weakly distributive and by using an earlier result of the author that makes it possible to eliminate the bound~$1$ in~$L$ (cf. Section~\ref{S:Conseq}), it follows that the semilattice $S$ is not isomorphic to $\Conc L$, for any lattice~$L$. Hence the ideal lattice of $S$ is not isomorphic to the congruence lattice of any lattice.

We shall now give some background on the problem solved by our theorem.
Funayama and Nakayama \cite{FuNa42} proved in 1942 that $\Con L$ is
\emph{distributive}, for any lattice~$\seq{L,\vee,\wedge}$. Dilworth proved soon after that
conversely, every finite distributive lattice is isomorphic to the
congruence lattice of some finite lattice (see \cite[pp.~455--456]{BFK} and \cite{GrDI}).
Birkhoff and Frink \cite{BiFr} proved in 1948 that the congruence lattice of any
algebra is what is nowadays called an \emph{algebraic} lattice, that is, it is
\emph{complete} and every element is a join of \emph{compact} elements
(see \cite{GLT2}). The question whether every algebraic distributive
lattice is isomorphic to~$\Con L$ for some lattice~$L$, often referred to
as CLP (`Congruence Lattice Problem'), is one of the most intriguing and
longest-standing open problems of lattice theory. In some sense, its first published occurrence is with the finite case as an exercise with asterisk (attributed to Dilworth) in the 1948 edition of Birkhoff's lattice theory book~\cite{Birk48}. The first published proof of this result seems to appear in Gr\"atzer and Schmidt's 1962 paper~\cite{GrSc62}. However, it seems that the earliest attempts at CLP were made by Dilworth himself, see \cite[pp.~455--456]{BFK}.

This problem has generated an enormous amount of work since then, in a
somewhat complex pattern of interconnected waves. Gr\"atzer and Schmidt
proved in 1963 that every algebraic lattice is isomorphic to the
congruence lattice of some algebra~\cite{GrSc63}. The reader can find in
Schmidt's monograph \cite{Schm82} a survey about congruence lattice
representations of algebras. The surveys by Gr\"atzer and Schmidt
\cite{GrScC,GrSc04} and Gr\"atzer's monograph \cite{FCBook} are focused on
congruence lattices of (mainly finite) lattices, while the survey by
T\r{u}ma and Wehrung \cite{CLPSurv} is more focused on congruence lattices
of infinite lattices. The main connection between the finite case and the
infinite case originates in Pudl\'ak's idea \cite{Pudl} of lifting, with
respect to the $\Conc$ functor, \emph{diagrams} of finite distributive
\jzs s. R\r{u}\v{z}i\v{c}ka, T\r{u}ma, and
Wehrung prove in \cite{RTW} that there are bounded lattices of cardinality
$\aleph_2$ whose congruence lattices are isomorphic neither to the normal
subgroup lattice of any group, nor to the submodule lattice of any module;
furthermore, the bound~$\aleph_2$ is optimal. Some of the more recent
works emphasize close connections between congruence lattices of lattices,
ideal lattices of rings, dimension theory of lattices, and nonstable
K-theory of rings, see for example
\cite{Berg86,GoWe1,GoWe2,Ruzi,NonMeas,WeURP,WReg}.

Distributive algebraic lattices are ideal lattices of distributive \jzs s
(see Section~\ref{S:Basic}), and for a lattice $L$, $\Con L$ is isomorphic to the ideal lattice of $\Conc L$. We obtain the following more convenient equivalent formulation of CLP (see \cite{CLPSurv} for details):

\begin{CLPS}
Is every distributive \jzs\ \emph{representable}, that is, isomorphic to
$\Conc L$, for some lattice $L$?
\end{CLPS}

In particular, the semilattice $S$ of our theorem provides a counterexample to~CLP.

Among the classical positive partial results are the following:
\begin{itemize}
\item[(1)] Every distributive \jzs\ $S$ of cardinality at most $\aleph_1$
is representable, see Huhn \cite{Huhn89a,Huhn89b}.

\item[(2)] Every distributive lattice with zero is representable, see
Schmidt \cite{Schm81}.
\end{itemize}

Further works extended the class of all representable distributive \jzs s,
for example to all \jz-direct limits of sequences of distributive
lattices with zero, see \cite{FPLat}. Moreover, the representing lattice
$L$ can be taken \emph{relatively complemented with zero}. This also
holds for case~(2) above. However, the latter result has been extended
further by R\r{u}\v{z}i\v{c}ka \cite{Ruzi}, who proved that the
representing lattice can be taken relatively complemented, modular, and
locally finite. This is not possible for~(1) above, as, for
$|S|\leq\aleph_1$, one can take $L$ relatively complemented
modular~\cite{WReg}, relatively complemented and locally
finite~\cite{GLW}, but not necessarily both~\cite{CXAl1}.

On the negative side, the works in \cite{PTW,TuWe,NonMeas,WeURP} show
that lattices with permutable congruences are not sufficient to solve
CLP. More precisely, \emph{there exists a representable distributive
\jzus\ of cardinality $\aleph_2$ that is not isomorphic to $\Conc L$ for
any lattice $L$ with permutable congruences}. The finite combinatorial
reason for this lies in the impossibility to prove certain `congruence
amalgamation properties'. The infinite combinatorial reason for this is
Kuratowski's Free Set Theorem (see Section~\ref{S:Basic}). The latter is
used to prove that certain infinitary statements called `uniform
refinement properties' fail in certain distributive semilattices.

Our proof carries a flavor of commutator theory with the structure of a semilattice, essentially because of Lemma~\ref{L:erosion}, the \emph{Erosion Lemma}. A precedent of this sort of situation occurs with Bill Lampe's wonderful
trick used in \cite{FLT} to prove that certain algebraic lattices
require, for their congruence representations, algebras with many
operations: namely, the term condition used in commutator theory in,
say, congruence-modular varieties (or larger, as considered in
\cite{KeSz98,Bowtie}).

\section{Basic concepts}\label{S:Basic}

A \jzs\ $S$ is \emph{distributive}, if $\bc\leq\ba\vee\bb$ in $S$ implies
that there are $\bx\leq\ba$ and $\by\leq\bb$ in $S$ such that
$\bx\leq\ba$, $\by\leq\bb$, and $\bc=\bx\vee\by$. Equivalently, the
ideal lattice of $S$ is a distributive lattice, see
\cite[Section~II.5]{GLT2}.

The assignment $L\mapsto\Conc L$ is extended, the usual way, to a functor from algebras with homomorphisms to \jzs s with \jzh s. For a positive integer~$m$, an algebra~$L$ has \emph{$(m+1)$-permutable congruences}, if $\ba\vee\bb=\bc_0\circ\bc_1\circ\cdots\circ\bc_m$ where $\bc_i$ equals~$\ba$ if~$i$ is even and~$\bb$ if~$i$ is odd, for all congruences~$\ba$ and~$\bb$ of~$L$ (the symbol $\circ$ denotes, as usual, composition of relations).

For an algebra $L$ endowed with a structure of
semilattice, with operation (thought of as a join operation) denoted by~$\vee$, we put $\Theta_L^+(x,y)= \Theta_L(y,x\vee y)$, for all $x,y\in L$. We say that the semilattice structure on $L$ is \emph{congruence-compatible}, if every congruence of~$L$ is a congruence for~$\vee$ (this definition extends to any operation instead of~$\vee$); equivalently, $x\equiv y\pmod{\ba}$ implies that $x\vee z\equiv y\vee z\pmod{\ba}$, for any $x,y,z\in L$ and any~$\ba\in\Con L$. In such a case, $\Theta_L^+(x,z)\subseteq\Theta_L^+(x,y)\vee\Theta_L^+(y,z)$, for any $x,y,z\in L$.

For partially ordered sets $P$ and $Q$, a map $f\colon P\to Q$ is
\emph{isotone}, if $x\leq y$ implies that $f(x)\leq f(y)$, for all
$x,y\in P$.

We shall also use standard set-theoretical notation and terminology,
referring the reader to \cite{Jech} for further information. We shall
denote by $\Pow(X)$ the powerset of a set $X$, by $[X]^{<\go}$ the set
of all finite subsets of $X$, and by $[X]^n$ (for $n<\go$) the set of
all $n$-element subsets of $X$. For a map
$\Phi\colon[X]^n\to[X]^{<\go}$, we say that an $(n+1)$-element subset
$U$ of $X$ is \emph{free with respect to $\Phi$}, if
$x\notin\Phi(U\setminus\set{x})$ for all $x\in U$. The following
statement of infinite combinatorics is one direction of a theorem due to
Kuratowski~\cite{Kura51}.

\begin{KFST}
Let $n$ be a natural number and let $X$ be a set with $|X|\geq\aleph_n$.
For every map $\Phi\colon[X]^n\to[X]^{<\go}$, there exists a
$(n+1)$-element free subset of $X$ with respect to $\Phi$.
\end{KFST}

We identify every natural number~$n$ with the set $\set{0,1,\dots,n-1}$,
and we denote by $\go$ the set of all natural numbers, which is also the
first limit ordinal. We shall usually denote elements in semilattices by
bold math characters $\ba$, $\bb$, $\bc$, \dots.

\section{Free distributive extension of a \jzs}\label{S:FDE}

As in \cite{NonPoset}, we shall use the construction of a ``free
distributive extension'' $\cR(S)$ of a \jzs\ $S$ given by Plo\v{s}\v{c}ica
and T\r{u}ma in \cite[Section~2]{PlTu}. The larger semilattice $\cR(S)$
is constructed by adding new elements $\bt(\ba,\bb,\bc)$, for
$\ba,\bb,\bc\in S$ such that $\bc\leq\ba\vee\bb$, subjected only to the
relations $\bc=\bt(\ba,\bb,\bc)\vee\bt(\bb,\ba,\bc)$ and
$\bt(\ba,\bb,\bc)\leq\ba$. It is a semilattice version of the
dimension group construction
$\mathbf{I}_K(E)$ presented in \cite[Section~1]{NonMeas}.
For convenience, we present an equivalent formulation here.

For a \jzs\ $S$, we shall put
$\cC(S)=\setm{\seq{\bu,\bv,\bw}\in S^3}{\bw\leq\bu\vee\bv}$.
A \emph{finite} subset $\bx$ of $\cC(S)$ is \emph{projectable} (resp.,
\emph{reduced}), if it satisfies condition~(1) (resp., (1)--(3)) below:
\begin{itemize}
\item[(1)] $\bx$ contains exactly one diagonal triple, that is, a triple
of the form $\seq{\bu,\bu,\bu}$; we put $\bu=\pi(\bx)$.

\item[(2)] $\seq{\bu,\bv,\bw}\in\bx$ and $\seq{\bv,\bu,\bw}\in\bx$ implies
that $\bu=\bv=\bw$, for all $\bu,\bv,\bw\in\nobreak S$.

\item[(3)]
$\seq{\bu,\bv,\bw}\in\bx\setminus\set{\seq{\pi(\bx),\pi(\bx),\pi(\bx)}}$
implies that $\bu,\bv,\bw\nleq\pi(\bx)$, for all $\bu,\bv,\bw\in S$.
\end{itemize}

In particular, observe that if $\bx$ is reduced,
$\seq{\bu,\bv,\bw}\in\bx$, and $\seq{\bu,\bv,\bw}$ is non-diagonal, then
$\bu\neq\bv$ and the elements $\bu$, $\bv$, and $\bw$ are nonzero.

We denote by $\ol{\cR}(S)$ (resp., $\cR(S)$) the set of all projectable
(resp., reduced) subsets of $\cC(S)$, endowed with the binary
relation $\leq$ defined by
 \begin{equation}\label{Eq:DefleqRed}
 \bx\leq\by\ \Longleftrightarrow\
 \forall\seq{\bu,\bv,\bw}\in\bx\setminus\by,
 \text{ either }\bu\leq\pi(\by)\text{ or }\bw\leq\pi(\by).
 \end{equation}
We call $\pi$ the \emph{canonical projection} from $\cR(S)$ onto $S$.
Observe that in general,~$\pi$ is not a \jh\ (however, see Remark~\ref{Rk:IdentifSRS}).
It is straightforward to verify that $\leq$ is a partial ordering on
$\ol{\cR}(S)$ (and thus on the subset $\cR(S)$). Now we shall present, in
terms of rewriting rules, the steps (i)--(iv) of the algorithm stated in \cite[Lemma~2.1]{PlTu}, aiming at Corollary~\ref{C:R(S)semil2}.

For finite subsets $\bx$ and $\by$ of $\cC(S)$, let $\bx\rightarrow_1\by$
hold, if there exists a non-diagonal $\seq{\ba,\bb,\bc}\in\bx$ such
that $\seq{\bb,\ba,\bc}\in\bx$ and
$\by=(\bx\setminus\set{\seq{\ba,\bb,\bc},\seq{\bb,\ba,\bc}})
\cup\set{\seq{\bc,\bc,\bc}}$. Denote by $\rightarrow_1^*$ the reflexive
and transitive closure of $\rightarrow_1$ on finite subsets of $\cC(S)$,
and denote by~$\cR_1(S)$ the set of all finite $\bx\subseteq\cC(S)$ such
that $\seq{\ba,\bb,\bc}\in\bx$ and $\seq{\bb,\ba,\bc}\in\bx$ implies that
$\ba=\bb=\bc$, for all $\ba,\bb,\bc\in S$. Put
$\ol{\cR}_1(S)=\ol{\cR}(S)\cap\cR_1(S)$. For a finite subset $\bx$ of
$\cC(S)$, we put
 \begin{multline*}
 \varphi(\bx)=(\bx\setminus\setm{\seq{\bu,\bu,\bu}}
 {\bu\in X})\cup\Set{\Seq{\bigvee X,\bigvee X,\bigvee X}},\\
 \text{where }X=\setm{\bu\in S}{\seq{\bu,\bu,\bu}\in\bx}.
 \end{multline*}
For $\bx\in\ol{\cR}(S)$ and a finite subset $\by$ of $\cC(S)$, let
$\bx\rightarrow_2\by$ hold, if there exists a non-diagonal
$\seq{\ba,\bb,\bc}\in\bx$ such that $\bb\leq\pi(\bx)$ and
 \[
 \by=(\bx\setminus\set{\seq{\ba,\bb,\bc},
 \seq{\pi(\bx),\pi(\bx),\pi(\bx)}})
 \cup\set{\seq{\bc\vee\pi(\bx),\bc\vee\pi(\bx),\bc\vee\pi(\bx)}}.
 \]
Observe that necessarily, $\by$ belongs to $\ol{\cR}(S)$ as well, and
denote by $\rightarrow_2^*$ the reflexive and transitive closure of
$\rightarrow_2$ on $\ol{\cR}(S)$. Denote by $\cR_2(S)$ the set of all
$\bx\in\ol{\cR}_1(S)$ such that for all non-diagonal
$\seq{\ba,\bb,\bc}\in\bx$, the inequality
$\bb\nleq\pi(\bx)$ holds. For any $\bx\in\ol{\cR}(S)$, we put
 \[
 \psi(\bx)=\bx\setminus
 \setm{\seq{\ba,\bb,\bc}\in\bx\text{ non-diagonal }}
 {\text{ either }\ba\leq\pi(\bx)\text{ or }\bc\leq\pi(\bx)}.
 \]
The correspondence with the algorithm stated in \cite[Lemma~2.1]{PlTu} is
as follows: the relation~$\rightarrow_1$ corresponds to step~(i); the
function~$\varphi$ corresponds to step~(ii); the relation~$\rightarrow_2$
corresponds to step~(iii); the function~$\psi$ corresponds to step~(iv).
The following lemma is a reformulation, in terms of $\rightarrow_1$,
$\rightarrow_2$, $\varphi$, and~$\psi$, of
\cite[Lemma~2.1]{PlTu}.

\begin{lemma}\label{L:R(S)semil1}
Let $\bx,\by\in\cR(S)$. Then there exists
$\seq{\bz_1,\bz_2}\in\cR_1(S)\times\cR_2(S)$ such that
$\bx\cup\by\rightarrow_1^*\bz_1$ and $\varphi(\bz_1)\rightarrow_2^*\bz_2$.
Furthermore, for any such pair $\seq{\bz_1,\bz_2}$, $\varphi(\bz_1)$
belongs to $\ol{\cR}_1(S)$ and $\psi(\bz_2)$ is the join,
in $\cR(S)$, of $\bx$ and $\by$.
\end{lemma}

\begin{corollary}\label{C:R(S)semil2}
The set $\cR(S)$ is a \jzs\ under the partial ordering defined in
\textup{\eqref{Eq:DefleqRed}}. Furthermore, the map
$j_S\colon S\to\cR(S)$, $\bx\mapsto\set{\seq{\bx,\bx,\bx}}$ is a \jze.
\end{corollary}

\begin{remark}\label{Rk:IdentifSRS}
We shall identify $\bx$ with the
element $\set{\seq{\bx,\bx,\bx}}$ of $\cR(S)$, for all $\bx\in\nobreak S$.
Then observe that the canonical map $\pi\colon\cR(S)\onto S$ is
\emph{isotone} and that the restriction of $\pi$ to $S$ is the identity.
The following is an easy consequence of \eqref{Eq:DefleqRed}.
 \begin{equation}\label{Eq:SrelcpleteR(S)}
 \bx\leq\by\ \Longleftrightarrow\ \bx\leq\pi(\by),
 \qquad\text{for all }\seq{\bx,\by}\in S\times\cR(S).
 \end{equation}
Now the elements of $\cR(S)\setminus S$ are exactly those
subsets~$\bx$ of~$\cC(S)\cup S$ (disjoint union) containing exactly one
element of~$S$, denoted by~$\pi(\bx)$, while
$\bx\setminus\set{\pi(\bx)}$ is nonempty and all its elements are triples
$\seq{\ba,\bb,\bc}\in\cC(S)$ such that $\seq{\bb,\ba,\bc}\notin\bx$ and
$\ba,\bb,\bc\nleq\pi(\bx)$.
\end{remark}

We shall use the symbol $\bt_S$, or $\bt$ if $S$ is understood, to denote
the elements of~$\cR(S)$ defined as
 \begin{equation*}
 \bt_S(\bu,\bv,\bw)=\begin{cases}
 \bw,&\text{if either }\bu=\bv\text{ or }\bv=0\text{ or }\bw=0,\\
 0,&\text{if }\bu=0,\\
 \set{\seq{0,0,0},\seq{\bu,\bv,\bw}},&\text{otherwise},
 \end{cases}
 \end{equation*}
for all $\seq{\bu,\bv,\bw}\in\cC(S)$. Then one can prove easily the
formula
 \begin{equation}\label{Eq:bxsupbabc}
 \bx=\bigvee\famm{\bt_S(\ba,\bb,\bc)}{\seq{\ba,\bb,\bc}\in\bx},
 \quad\text{for all }\bx\in\cR(S).
 \end{equation}
The following is a slight strengthening of \cite[Theorem~2.3]{PlTu}, with
the same proof. The uniqueness statement follows from
\eqref{Eq:bxsupbabc}.

\begin{lemma}\label{L:Interp}
Let $S$ and $T$ be \jzs s and let $f\colon S\to T$ be a \jzh.
Furthermore, let $\imath\colon\cC(\im f)\to T$ be a map such that
$\imath(\bx,\by,\bz)\vee\imath(\by,\bx,\bz)=\bz$ and
$\imath(\bx,\by,\bz)\leq\bx$, for all $\seq{\bx,\by,\bz}\in\cC(\im f)$.
Then there exists a unique map $f_{(\imath)}\colon\cR(S)\to T$ such that
$f_{(\imath)}(\bt_S(\bx,\by,\bz))=\imath(f(\bx),f(\by),f(\by))$, for all
$\seq{\bx,\by,\bz}\in\cC(S)$.
\end{lemma}

By applying Lemma~\ref{L:Interp} to the map $j_T\circ f$ and
defining $\imath$ as the restriction of~$\bt_T$ to~$\cC(\im f)$,
we obtain item~(1) of the following result. Item~(2) follows easily.
\goodbreak

\begin{proposition}\label{P:StR(S)}\hfill
\begin{enumerate}
\item For \jzs s $S$ and $T$, every \jzh\ $f\colon S\to T$ extends to a
unique \jzh\ $\cR(f)\colon\cR(S)\to\cR(T)$ such that
$\cR(f)(\bt_S(\bu,\bv,\bw))=\bt_T(f(\bu),f(\bv),f(\bw))$, for all
$\seq{\bu,\bv,\bw}\in\cC(S)$.

\item The assignment $S\mapsto\cR(S)$, $f\mapsto\cR(f)$ is a functor.
\end{enumerate}
\end{proposition}

Putting $\cR^0(S)=S$ and
$\cR^{n+1}(S)=\cR(\cR^n(S))$ for each~$n$, the increasing
union $\cD(S)=\bigcup\famm{\cR^n(S)}{n<\go}$ is a distributive \jzs\
extending~$S$. Furthermore,
putting $\cD(f)=\bigcup\famm{\cR^n(f)}{n<\go}$ for each \jzh\ $f$, we
obtain that~$\cD$ is a functor.
The proof of the following lemma is straightforward.

\begin{lemma}\label{L:RD(inters)}
Let $S$ be a \jzs\ and let $\seqm{S_i}{i\in I}$ be a family of
\jz-subsemilattices of $S$. The following statements hold:
\begin{enumerate}
\item $\cR\left(\bigcap_{i\in I}S_i\right)=\bigcap_{i\in I}\cR(S_i)$
and $\cD\left(\bigcap_{i\in I}S_i\right)=\bigcap_{i\in I}\cD(S_i)$.
\item If $I$ is a nonempty upward directed partially ordered set and
$\seqm{S_i}{i\in I}$ is isotone, then
$\cR\left(\bigcup_{i\in I}S_i\right)=\bigcup_{i\in I}\cR(S_i)$ and
$\cD\left(\bigcup_{i\in I}S_i\right)=\bigcup_{i\in I}\cD(S_i)$.
\end{enumerate}
\end{lemma}

\begin{definition}\label{D:Cplxx}
For a \jzs\ $S$ and an element $\bx\in\cD(S)$, we define the \emph{rank}
of $\bx$, denoted by $\rk\bx$, as the least natural number $n$ such that
$\bx\in\cR^n(S)$.
\end{definition}

\section{The functors $\cL$ and $\cG$}\label{S:cLcG}

In the present section we shall construct the semilattice used in the counterexample and demonstrate one of its crucial properties, namely the \emph{Evaporation Lemma} (Lemma~\ref{L:Evap}).

For a set $\Omega$, we denote by $\cL(\Omega)$ the \jzs\ defined by
generators $1$ and $\ba_0^{\xi}$, $\ba_1^{\xi}$ (for $\xi\in\Omega$),
subjected to the relations
 \begin{equation}\label{Eq:Rela0a1xi}
 \ba_0^{\xi}\vee\ba_1^{\xi}=1,\quad\text{for all }\xi\in\Omega.
 \end{equation}
Hence $\cL(\Omega)$ is the same semilattice as the one presented in
\cite[Section~3]{PlTu}. It is a semilattice version of the dimension
group $\mathbf{E}_K(\Omega)$ presented in \cite[Section~2]{NonMeas}.
It can be `concretely'
represented as the (semi)lattice of all pairs
$\seq{X,Y}\in\Pow(\Omega)\times\Pow(\Omega)$ such that either $X$ and $Y$
are finite and disjoint or $X=Y=\Omega$, with
 \[
 \ba_0^{\xi}=\seq{\set{\xi},\es}\text{ and }
 \ba_1^{\xi}=\seq{\es,\set{\xi}},\quad\text{for all }\xi\in\Omega.
 \]
We shall identify $\cL(X)$ with the \jzu-subsemilattice of $\cL(\Omega)$
generated by the subset
$\setm{\ba_i^{\xi}}{\xi\in X\text{ and }i<2}$, for all
$X\subseteq\Omega$. For sets $X$ and $Y$, any map $f\colon X\to Y$ gives
rise to a unique \jzuh\ $\cL(f)\colon\cL(X)\to\cL(Y)$ such that
$\cL(f)(\ba_i^{\xi})=\ba_i^{f(\xi)}$, for all
$\seq{\xi,i}\in X\times\set{0,1}$. Of course, the assignment
$X\mapsto\cL(X)$, $f\mapsto\cL(f)$ is a functor from the category of sets
with maps to the category of \jzus s and
\jzuh s.

Next, we put $\cG=\cD\circ\cL$, the composition of the two functors $\cD$
and $\cL$. Hence, for a set $\Omega$, the semilattice $\cG(\Omega)$ may
be loosely described as a `free distributive \jzs\ defined by generators
$\ba_i^{\xi}$, for $\xi\in\Omega$ and $i<2$, and relations
\eqref{Eq:Rela0a1xi}'. It is a distributive \jzus, of the same cardinality
as $\Omega$ in case $\Omega$ is infinite.

The proof of the following lemma is straightforward (see
Lemma~\ref{L:RD(inters)}).

\begin{lemma}\label{L:CapCupLGpres}
Let $\Omega$ be a set and let $\seqm{X_i}{i\in I}$ be a family of subsets
of $\Omega$. The following statements hold:
\begin{enumerate}
\item $\cL\bigl(\bigcap_{i\in I}X_i\bigr)=\bigcap_{i\in I}\cL(X_i)$ and
$\cG\bigl(\bigcap_{i\in I}X_i\bigr)=\bigcap_{i\in I}\cG(X_i)$.

\item If $I$ is a nonempty upward directed partially ordered set and the family
$\seqm{X_i}{i\in I}$ is isotone, then
$\cL\bigl(\bigcup_{i\in I}X_i\bigr)=\bigcup_{i\in I}\cL(X_i)$ and
$\cG\bigl(\bigcup_{i\in I}X_i\bigr)=\bigcup_{i\in I}\cG(X_i)$.
\end{enumerate}
\end{lemma}

\begin{corollary}\label{C:ExistSupp}
For any set $\Omega$ and any $\bx\in\cG(\Omega)$, there exists a
least \pup{finite} subset~$X$ of~$\Omega$ such that $\bx\in\cG(X)$.
\end{corollary}

We shall call the subset $X$ of Corollary~\ref{C:ExistSupp} the
\emph{support} of $\bx$, and denote it by $\supp(\bx)$.

\begin{lemma}\label{L:xleqy+ai}
Let $\Omega$ be a set, let $\alpha\in\Omega$, and let $i<2$. Then
$\bx\leq\by\vee\ba_i^{\alpha}$ implies that $\bx\leq\by$, for all
$\bx,\by\in\cG(\Omega\setminus\set{\alpha})$.
\end{lemma}

\begin{proof}
There exists a unique retraction
$r\colon\cL(\Omega)\onto\cL(\Omega\setminus\set{\alpha})$ such that
$r(\ba_i^{\alpha})=0$. Put $s=\cD(r)$, and observe that $s(\bx)=\bx$,
$s(\by)=\by$, and $s(\ba_i^{\alpha})=0$. By applying $s$ to the inequality
$\bx\leq\by\vee\ba_i^{\alpha}$, we get the conclusion.
\end{proof}

The following crucial lemma describes an `evaporation process' in
$\cG(\Omega)$.

\begin{lemma}[Evaporation Lemma]\label{L:Evap}
Let $\alpha$, $\beta$, $\delta$ be distinct elements in a set
$\Omega$, let $i,j<2$, $\bx\in\cG(\Omega\setminus\set{\beta})$,
$\by\in\cG(\Omega\setminus\set{\alpha})$, and
$\bz\in\cG(\Omega\setminus\set{\delta})$. Then 
 \[
 \bz\leq\bx\vee\by,\quad\bx\leq\ba_0^{\delta},\ba_i^{\alpha},
 \quad\text{and}\quad
 \by\leq\ba_1^{\delta},\ba_j^{\beta}
 \]
implies that $\bz=0$.
\end{lemma}

\begin{proof}
For $s\in\go$ and $\bu\in\cR^{s+1}\cL(\Omega)\setminus\cR^s\cL(\Omega)$,
we shall denote by $\pi(\bu)$ the image of $\bu$ under the canonical projection from
$\cR^{s+1}\cL(\Omega)$ to $\cR^s\cL(\Omega)$.
Put $m=\rk\bx$, $n=\rk\by$, and $k=\rk\bz$. We argue by induction on
$m+n+k$. If $\bz\leq\bx$, then $\bz\leq\ba_0^{\delta}$, thus, as
$\bz\in\cG(\Omega\setminus\set{\delta})$, it follows from
Lemma~\ref{L:xleqy+ai} that $\bz=0$ so we are done. The conclusion is
similar in case $\bz\leq\by$. So suppose that $\bz\nleq\bx,\by$. If
$m=0$, then, as $\bx\in\cL(\Omega)$ and
$\bx\leq\ba_0^{\delta},\ba_i^{\alpha}$ with $\alpha\neq\delta$, we get
$\bx=0$, so $\bz\leq\by$, \contr; hence $m>0$. Similarly, $n>0$. Put
$l=\max\set{m,n}$, $\bx^*=\bx\setminus\set{\pi(\bx)}$, and
$\by^*=\by\setminus\set{\pi(\by)}$ (see Remark~\ref{Rk:IdentifSRS}).
Furthermore, we define (using again Remark~\ref{Rk:IdentifSRS})
a finite subset $\bw$ of $\cC\cR^{l-1}\cL(\Omega)$ as
 \begin{equation}\label{Eq:Defw}
 \bw=\begin{cases}
 \bx^*\cup\by^*\cup\set{\pi(\bx)\vee\pi(\by)},&\text{if }m=n,\\
 \by^*\cup\set{\bx\vee\pi(\by)},&\text{if }m<n,\\
 \bx^*\cup\set{\pi(\bx)\vee\by},&\text{if }m>n.
 \end{cases}
 \end{equation}

\begin{sclaim}
The set $\bw$ belongs to $\cR^l\cL(\Omega)$, and $\bx,\by\leq\bw$.
\end{sclaim}

\begin{scproof}
We need to verify that $\bw$ is a reduced subset of
$\cC\cR^{l-1}\cL(\Omega)$, modulo the identification of elements with
diagonal triples (see Remark~\ref{Rk:IdentifSRS}).
It is obvious that there exists exactly one element in
$\bw\cap\cR^{l-1}\cL(\Omega)$, namely,
 \[
 \pi(\bw)=\begin{cases}
 \pi(\bx)\vee\pi(\by),&\text{if }m=n,\\
 \bx\vee\pi(\by),&\text{if }m<n,\\
 \pi(\bx)\vee\by,&\text{if }m>n.
 \end{cases}
 \]
This settles item~(1) of the definition of a reduced set.

Now suppose that there exists a non-diagonal triple $\seq{\ba,\bb,\bc}$
of elements of $\cR^{l-1}\cL(\Omega)$ such that $\seq{\ba,\bb,\bc}\in\bw$
and $\seq{\bb,\ba,\bc}\in\bw$. As both $\bx$ and $\by$ are reduced sets,
the only possibility is $m=n$ and, say, $\seq{\ba,\bb,\bc}\in\bx$
and $\seq{\bb,\ba,\bc}\in\by$. As $\bx\in\cG(\Omega\setminus\set{\beta})$
and $\by\in\cG(\Omega\setminus\set{\alpha})$, all elements $\ba$,
$\bb$, $\bc$ belong to $\cG(\Omega\setminus\set{\alpha,\beta})$
(see Lemma~\ref{L:CapCupLGpres}).
As $\seq{\ba,\bb,\bc}\in\bx$ and $\bx\leq\ba_i^{\alpha}$, it follows
from \eqref{Eq:DefleqRed} and the assumption that $\seq{\ba,\bb,\bc}$ is
non-diagonal that either $\ba\leq\ba_i^{\alpha}$ or
$\bc\leq\ba_i^{\alpha}$. As $\ba,\bc\in\cG(\Omega\setminus\set{\alpha})$,
it follows from Lemma~\ref{L:xleqy+ai} that either $\ba=0$ or $\bc=0$,
\contr. This settles item~(2) of the definition of a reduced set.

Finally, let $\seq{\ba,\bb,\bc}\in\bw$ be a non-diagonal triple of
elements of $\cR^{l-1}\cL(\Omega)$, we must verify that
$\ba,\bb,\bc\nleq\pi(\bw)$. Suppose, for example, that $\ba\leq\pi(\bw)$.
If $m=n$, then $\ba\leq\pi(\bx)\vee\pi(\by)$ and, say,
$\seq{\ba,\bb,\bc}\in\bx^*$. {}From $\pi(\by)\leq\by\leq\ba_j^{\beta}$ it
follows that $\ba\leq\pi(\bx)\vee\ba_j^{\beta}$. As
$\ba,\pi(\bx)\in\cG(\Omega\setminus\set{\beta})$ and by
Lemma~\ref{L:xleqy+ai}, it follows that $\ba\leq\pi(\bx)$, which
contradicts the assumption that $\seq{\ba,\bb,\bc}$ is a non-diagonal
triple in~$\bx$. If $m<n$, then $\seq{\ba,\bb,\bc}\in\by^*$ and
$\ba\leq\bx\vee\pi(\by)$, so $\ba\leq\ba_i^{\alpha}\vee\pi(\by)$, and so,
as $\ba,\pi(\by)\in\cG(\Omega\setminus\set{\alpha})$ and by
Lemma~\ref{L:xleqy+ai}, it follows that
$\ba\leq\pi(\by)$, which
contradicts the assumption that $\seq{\ba,\bb,\bc}$ is a non-diagonal
triple in~$\by$. The proof for the case $m>n$ is similar. So we have
proved that $\ba\nleq\pi(\bw)$. The proofs for~$\bb$ and~$\bc$ are
similar. This settles item~(3) of the definition of a reduced set.

The verification of the inequalities $\bx,\by\leq\bw$ (see
\eqref{Eq:DefleqRed}) is straightforward.
In fact, it is not hard to verify, using Lemma~\ref{L:R(S)semil1}, that
$\bw=\bx\vee\by$.
\end{scproof}

Now we complete the proof of Lemma~\ref{L:Evap}. {}From the claim above
it follows that $\bz\leq\bw$. If $k<l$ then $\bz\leq\pi(\bw)$, hence, as
$\pi(\bw)\in\set{\pi(\bx)\vee\pi(\by),\bx\vee\pi(\by),\pi(\bx)\vee\by}$
and by the induction hypothesis, $\bz=0$. So suppose from now on that
$k\geq l$; in particular, $k>0$. As $\pi(\bz)\leq\bz\leq\bx\vee\by$, it
follows from the induction hypothesis that $\pi(\bz)=0$. Hence, if
$\bz\neq0$, then there exists a non-diagonal triple
$\seq{\ba,\bb,\bc}\in\bz\cap\cC\cR^{l-1}\cL(\Omega)$. As $\bz\leq\bw$, we
obtain that either
$\seq{\ba,\bb,\bc}\in\bw$ or $\ba\leq\bw$ or $\bc\leq\bw$. In the first
case, say, $\seq{\ba,\bb,\bc}\in\bx$, we get
$\bt(\ba,\bb,\bc)\leq\bx\leq\ba_0^{\delta}$ with
$\ba,\bb,\bc\in\cG(\Omega\setminus\set{\delta})$ (because
$\seq{\ba,\bb,\bc}\in\bz$), so $\bt(\ba,\bb,\bc)=0$ by
Lemma~\ref{L:xleqy+ai}, \contr. If either $\ba\leq\bw$ or $\bc\leq\bw$,
then, by the induction hypothesis, either $\ba=0$ or $\bc=0$, \contr.
Therefore, $\bz=0$.
\end{proof}

\section{The Erosion Lemma}\label{S:EvapEros}

The proofs of our negative results are based on the conflict between a
\emph{non-structure theorem} on the semilattices $\cG(\Omega)$, here the `Evaporation Lemma' (Lem\-ma~\ref{L:Evap}), and a \emph{structure theorem} on arbitrary bounded semilattices, Lemma~\ref{L:erosion}, that we shall now introduce. This lemma, the \emph{Erosion Lemma}, contains, despite its extreme simplicity, the gist of the present paper. Moreover, further extensions of our methods seem to use the same formulation of the Erosion Lemma, while there seem to be many different `Evaporation Lemmas' (such as Lemma~\ref{L:Evap}).

{}From now on, we shall denote by $\eps$ the `parity function' on the natural numbers, defined by the rule
 \begin{equation}\label{Eq:Defeps}
 \eps(n)=\begin{cases}
 0,&\text{if }n\text{ is even},\\
 1,&\text{if }n\text{ is odd},
  \end{cases}
  \qquad\text{for every natural number }n.
 \end{equation}
Throughout this section, we let $L$ be an algebra possessing a congruence-compatible structure of semilattice $\seq{L,\vee}$. We put
 \[
 U\vee V=\setm{u\vee v}{\seq{u,v}\in U\times V},\quad
 \text{for all }U,V\subseteq L,
 \]
and we denote by $\Concs^UL$ the \jz-subsemilattice of $\Conc L$ generated by all principal congruences $\Theta_L(u,v)$, where $\seq{u,v}\in U\times U$.

\begin{lemma}[The Erosion Lemma]\label{L:erosion}
Let $x_0,x_1\in L$, and let $Z=\setm{z_i}{0\leq i\leq n}$, with $n\in\gos$, be a finite subset of~$L$ with $\bigvee_{i<n}z_i\leq z_n$. Put
 \[
 \ba_j=\bigvee\famm{\Theta_L(z_i,z_{i+1})}{i<n,\ \eps(i)=j},\text{ for all }j<2.
 \]
Then there are congruences $\bu_j\in\Concs^{\set{x_j}\vee Z}L$, for $j<2$, such that
 \[
 z_0\vee x_0\vee x_1\equiv z_n\vee x_0\vee x_1
 \pmod{\bu_0\vee\bu_1}\quad\text{and}\quad
 \bu_j\subseteq\ba_j\cap\Theta_L^+(z_n,x_j),\text{ for all }j<2.
 \]
\end{lemma}

\begin{proof}
Put $\bv_i=\Theta_L(z_i\vee x_{\eps(i)},z_{i+1}\vee x_{\eps(i)})$, for all $i<n$.
Observe that~$\bv_i$ belongs to $\Concs^{\set{x_{\eps(i)}}\cup Z}L$.
{}From $z_n\leq x_{\eps(i)}\pmod{\Theta_L^+(z_n,x_{\eps(i)})}$ and
$z_i\equiv z_{i+1}\pmod{\ba_{\eps(i)}}$ it follows, respectively
(and using $z_i\vee z_n=z_{i+1}\vee z_n$ in the first case), that
 \begin{equation}\label{Eq:bviThetaL}
 \bv_i\subseteq\Theta_L^+(z_n,x_{\eps(i)})\quad\text{and}\quad
 \bv_i\subseteq\ba_{\eps(i)}.
 \end{equation}
Now we put
 \[
 \bu_j=\bigvee\famm{\bv_i}{i<n,\ \eps(i)=j},\quad\text{for all }j<2.
 \]
Hence $\bu_j\in\Concs^{\set{x_j}\vee Z}L$, for all $j<2$. Furthermore,
from \eqref{Eq:bviThetaL} it follows that
$\bu_j\subseteq\ba_j\cap\Theta_L^+(z_n,x_j)$. Finally, from
$z_i\vee x_{\eps(i)}\equiv z_{i+1}\vee x_{\eps(i)}\pmod{\bv_i}$, for all $i<n$, it follows that
$z_i\vee x_0\vee x_1\equiv z_{i+1}\vee x_0\vee x_1\pmod{\bu_0\vee\bu_1}$. Therefore,
$z_0\vee x_0\vee x_1\equiv z_n\vee x_0\vee x_1\pmod{\bu_0\vee\bu_1}$.
\end{proof}

\section{The proof}\label{S:Proof}

Our main theorem is the following.

\begin{theorem}\label{T:NonRepr}
Let $\Omega$ be a set of cardinality at least $\aleph_{\omega+1}$ and
let $L$ be an algebra. If~$L$ has a congruence-compatible structure of
\jus, then there is no weakly distributive \jzh\ from $\Conc L$ to
$\cG(\Omega)$ with~$1$ in its range.
\end{theorem}

The remainder of this section will be devoted to a proof of
Theorem~\ref{T:NonRepr}. Suppose, to the contrary, that $L$ and
$\mu\colon\Conc L\to\cG(\Omega)$ are as above. We fix a congruence-compatible structure of \jus\ on~$L$. There are a positive integer~$m$ and elements $t_0$, \dots, $t_{m-1}$ in~$L$ such that
 \begin{equation}\label{Eq:Decomp1mu}
 \bigvee_{r<m}\mu\Theta_L(t_r,1)=1.
 \end{equation}
For each $\xi\in\Omega$, as $\mu\Theta_L(t_r,1)\leq1=\ba_0^{\xi}\vee\ba_1^{\xi}$ holds for each~$r<m$, we obtain, by using the weak distributivity of~$\mu$
at~$\Theta_L(t_r,1)$, an integer~$n_\xi\geq2$ and elements $z_{r,i}^{\xi}\in L$,
for $0\leq r<m$ and $0\leq i\leq n_{\xi}$, such that
$z_{r,0}^{\xi}=t_r$, $z_{r,n_{\xi}}^{\xi}=1$, and
 \begin{equation}\label{Eq:CongSmall}
 \mu\Theta_L(z_{r,i}^{\xi},z_{r,i+1}^{\xi})\leq\ba_{\eps(i)}^{\xi},
 \quad\text{for all }r<m\text{ and }i<n_{\xi}.
 \end{equation}
(We recall that $\eps$ is the parity function defined in \eqref{Eq:Defeps}.)
After replacing $z_{r,i}^\xi$ by $t_r\vee z_{r,i}^\xi$, we may also assume that $t_r\leq z_{r,i}^\xi$ holds, for all $r<m$, $i\leq n_\xi$, and $\xi\in\Omega$.
As $|\Omega|\geq\aleph_{\omega+1}$ and $\aleph_{\omega+1}$ is a regular
cardinal (this is the reason why $\aleph_{\omega}$ would not work
\emph{a priori}), there are a positive integer~$n$ and $\Omega'\subseteq\Omega$ such that $|\Omega'|=\aleph_{\omega+1}$ and $n_{\xi}=n$ for all $\xi\in\Omega'$.
Pick any retraction $\rho\colon\Omega\onto\Omega'$ and replace
$\mu$ by $\cG(\rho)\circ\mu$. We might lose the weak distributivity of $\mu$,
but we keep the elements $z_{r,i}^{\xi}$ and the statements \eqref{Eq:CongSmall}, which are all that matters. Furthermore, after replacing~$L$ by~$L/{\theta}$ where
$\seq{x,y}\in\theta$ if{f} $\mu\Theta_L(x,y)=0$ (for all $x,y\in L$), we
may assume that~\emph{$\mu$ separates zero}, that is,
$\mu^{-1}\set{0}=\set{0}$.

Hence we shall assume, from now on, that $\mu$ separates zero and
$n_{\xi}=n$ for all $\xi\in\Omega$. For every finite subset $X$ of
$\Omega$, we shall denote by $S(X)$ the join-subsemilattice of $L$
generated by $\setm{z_{r,i}^{\xi}}{0\leq r<m,\ 0\leq i\leq n,\text{ and }\xi\in X}$.
As $S(X)$ is finite, $\Phi(X)=\bigcup\Famm{\supp\mu\Theta_L(x,y)}{x,y\in S(X)}$ is a finite subset of $\Omega$.

As $|\Omega|\geq\aleph_{2^n}$, it follows from Kuratowski's Free Set Theorem that there exists a $(2^n+1)$-element subset $U$ of $\Omega$ which is free with respect to the restriction of~$\Phi$ to $2^n$-elements subsets of~$\Omega$.

For all natural numbers $k$, $l$ with $k\leq n-1$ and $l\leq 2^k$,
let $P(k,l)$ hold, if for all $r<m$ and all disjoint $X,Y\subseteq U$ with $|X|=2^k-l$
and $|Y|=2l$, the following equality $E_r(X,Y)$ holds:
 \begin{equation}
 \bigvee\famm{z_{r,n-k}^{\xi}}{\xi\in X}\vee
 \bigvee\famm{z_{r,n-k-1}^{\eta}}{\eta\in Y}=1.\tag{$E_r(X,Y)$}
 \end{equation}

The method used to prove Lemma~\ref{L:Descent} below could be described as `the erosion method': namely, prove, using the Erosion Lemma, that joins of larger and larger subsets of $L$ of the form
$\setm{z_{r,n-k}^{\xi}}{\xi\in X}\cup\setm{z_{r,n-k-1}^{\eta}}{\eta\in Y}$, with $k$ larger and larger, remain equal to~$1$. For large enough $k$, this will lead naturally to $t_r=1$.

\begin{lemma}[Descent Lemma]\label{L:Descent}
The statement $P(k,l)$ holds, for all natural numbers~$k,l$ such that $k\leq n-1$ and $l\leq 2^k$.
\end{lemma}

\begin{proof}
We argue by induction on $2^k+l$. Obviously, $P(0,0)$ holds.
Assuming that $P(k,l)$ holds, we shall establish $P(k',l')$ for the next
value $\seq{k',l'}$. As $P(k,2^k)$ is equivalent to
$P(k+1,0)$, we may assume that $l<2^k$, so $k'=k$ and $l'=l+1$. So
let $X,Y\subseteq U$ disjoint with $|X|=2^k-l-1$ and $|Y|=2l+2$. As
$|X|+|Y|=2^k+l+1\leq 2^n$ and $|U|=2^n+1$, there exists an element
$\delta\in U\setminus(X\cup Y)$. Pick $r<m$ and distinct elements
$\eta_0,\eta_1\in Y$, set $Y'=Y\setminus\set{\eta_0,\eta_1}$ and
 \begin{equation}\label{Eq:Defxj}
 x_j=\bigvee\famm{z_{r,n-k}^{\xi}}{\xi\in X}\vee
 \bigvee\famm{z_{r,n-k-1}^{\eta}}{\eta\in Y'\cup\set{\eta_j}},\quad\text{for all }j<2.
 \end{equation}
It follows from the induction hypothesis that
 \begin{equation}\label{Eq:nearxj=1}
 \bigvee\famm{z_{r,n-k}^{\xi}}{\xi\in X\cup\set{\eta_j}}\vee
 \bigvee\famm{z_{r,n-k-1}^{\eta}}{\eta\in Y'}=1,\quad\text{for all }j<2.
 \end{equation}
Now recall that, by \eqref{Eq:CongSmall},
 \[
 \mu\Theta_L(z_{r,n-k}^{\eta_j},z_{r,n-k-1}^{\eta_j})\leq\ba_{\eps(n-k-1)}^{\eta_j},
 \quad\text{for all }j<2.
 \]
Using \eqref{Eq:Defxj} and \eqref{Eq:nearxj=1}, it follows that
$\mu\Theta_L(x_j,1)\leq\ba_{\eps(n-k-1)}^{\eta_j}$, for all $j<2$. Therefore, using
Lemma~\ref{L:erosion} with $z_{r,i}^{\delta}$ in place of $z_i$,
for $0\leq i\leq n$, and observing that $t_r\leq x_0\vee x_1$
(because $t_r\leq z_{r,i}^\xi$ everywhere), we obtain congruences
$\bu_j\in\Concs^{S(X\cup Y'\cup\set{\eta_j,\delta})}L$, for $j<2$, such that
 \begin{equation}\label{Eq:muThx0x1small}
 \Theta_L(x_0\vee x_1,1)\leq\bu_0\vee\bu_1\quad\text{and}\quad
 \mu(\bu_j)\leq\ba_{\eps(n-k-1)}^{\eta_j},\ba_j^{\delta},\quad\text{for all }j<2.
 \end{equation}
It follows from the definition of $\Phi$ that
$\mu(\bu_j)\in\cG\Phi(X\cup Y'\cup\set{\eta_j,\delta})$ and\linebreak
$\mu\Theta_L(x_0\vee x_1,1)\in\cG\Phi(X\cup Y)$.
Using the monotonicity of $\Phi$ and the freeness of $U$ with respect to
the restriction of $\Phi$ to $2^n$-element subsets, we obtain
 \begin{align*}
 \Phi(X\cup Y)&\subseteq\Omega\setminus\set{\delta},\\
 \Phi(X\cup Y'\cup\set{\eta_j,\delta})&\subseteq
 \Omega\setminus\set{\eta_{1-j}},\quad\text{for all }j<2.
 \end{align*}
As $\mu\Theta_L(x_0\vee x_1,1)$ belongs to $\cG\Phi(X\cup Y)$ and by using
\eqref{Eq:muThx0x1small} together with Lemma~\ref{L:Evap}, we obtain that
$\mu\Theta_L(x_0\vee x_1,1)=0$, that is, since $\mu$ separates zero,
$x_0\vee x_1=1$, which completes the proof of the equality $E_r(X,Y)$.
\end{proof}

Now pick $\delta\in U$ and put $Y=U\setminus\set{\delta}$, so
$|Y|=2^n$. By applying Lemma~\ref{L:Descent} to $k=n-1$ and
$l=2^{n-1}$, we obtain the equality $\bigvee\famm{z_{r,0}^{\eta}}{\eta\in Y}=1$, that is, $t_r=1$. But this holds for all $r<m$, which contradicts~\eqref{Eq:Decomp1mu}. This
completes the proof of Theorem~\ref{T:NonRepr}.

\begin{remark}\label{Rk:WDattr}
In the assumptions of Theorem~\ref{T:NonRepr}, it is sufficient to restrict the weak distributivity assumption of~$\mu$ to congruences~$\Theta_L(t_r,1)$, for $r<m$, satisfying \eqref{Eq:Decomp1mu}.
\end{remark}

\section{Consequences on congruence lattices of lattices}\label{S:Conseq}

Observe that Theorem~\ref{T:NonRepr} applies to $L$ a \emph{lattice with a largest element}. We now extend this result to arbitrary lattices.

\begin{theorem}\label{T:NowdIm}
For any set $\Omega$ and any algebra~$L$ with a congruence-compatible lattice structure, if $|\Omega|\geq\aleph_{\go+1}$, then there exists no weakly distributive \jzh\ $\mu\colon\Conc L\to\cG(\Omega)$ with $1$ in its range.\end{theorem}

\begin{proof}
Denote by $L^{\mathrm{lat}}$ the given congruence-compatible lattice structure on (the underlying set of)~$L$. It is straightforward to verify that the canonical homomorphism from $\Conc(L^{\mathrm{lat}})$ to~$\Conc L$, that to each compact congruence of~$L^{\mathrm{lat}}$ associates the congruence of~$L$ that it generates, is weakly distributive. As the composition of two weakly distributive homomorphisms is weakly distributive, it suffices to prove the theorem in case~$L$ is a lattice.

So let $\mu\colon\Conc L\to\cG(\Omega)$ be a weakly distributive
\jzh\ with~$1$ in its range, where $|\Omega|\geq\aleph_{\go+1}$. As
$1=\bigvee_{i<n}\mu\Theta_L(u_i,v_i)$, for a positive integer~$n$ and elements $u_i\leq v_i$ in~$L$, for $i<n$, we get $1=\mu\Theta_L(u,v)$, where $u=\bigwedge_{i<n}u_i$ and $v=\bigvee_{i<n}v_i$. Put $K=[u,v]$. It follows from \cite[Proposition~1.2]{WeURP} that the canonical homomorphism $\jmath\colon\Conc K\to\Conc L$ is weakly distributive. Hence $\mu\circ\jmath$ is a weakly distributive homomorphism from~$\Conc K$ to~$\cG(\Omega)$ with~$1$ in its range, with~$K$ a \emph{bounded} lattice. This contradicts Theorem~\ref{T:NonRepr}.
\end{proof}

In particular, we obtain a negative solution to CLP.

\begin{corollary}\label{C:notCLP}
Let $\Omega$ be a set. If $|\Omega|\geq\aleph_{\go+1}$, then there exists no lattice $L$ with $\Conc L\cong\cG(\Omega)$.
\end{corollary}

By contrast, Lampe proved in \cite{Lamp82}
that \emph{every \jzus\ is isomorphic to $\Conc G$ for some groupoid $G$
with $4$-permutable congruences}. In particular, $\cG(\aleph_{\go+1})\cong\Conc G$ for some groupoid $G$ with $4$-permutable congruences, while there is no lattice $L$ such that $\cG(\aleph_{\go+1})\cong\Conc L$.
This shows a critical discrepancy between general algebras and
lattices.

\section{Discussion}\label{S:Disc}

\subsection{A new uniform refinement property}\label{Su:URP}
In many works such as \cite{PTW,RTW,TuWe,CLPSurv,WeURP,CXAl1,NonPoset}, the classes of semilattices that are representable with respect to various functors are separated from the corresponding counterexamples by infinitary statements called \emph{uniform refinement properties}. We shall now discuss briefly how this can also be done here. As the proofs do not seem to add much to the already existing results, we shall omit the details.

For a positive integer $m$ and a nonempty set $\Omega$, denote by $\Sem(m,\Omega)$ the join-semilattice defined by generators $\zero$, $\one$, and $k\cdot\dxi$ for
$0\leq k\leq m+1$ and $\xi\in\Omega$, subjected to the relations
 \[
 \zero=0\cdot\dxi\leq 1\cdot\dxi\leq\cdots\leq m\cdot\dxi\leq(m+1)\cdot\dxi=\one,
 \quad\text{for }\xi\in\Omega.
 \]

\begin{definition}\label{D:CLR}
For an element $\be$ in a \jzs\ $S$, we say that $S$ satisfies
$\CLR(\be)$, if for every nonempty set $\Omega$ and every family
$\seqm{\ba_i^{\xi}}{\seq{\xi,i}\in\Omega\times\set{0,1}}$ with entries in~$S$ such that $\be\leq\ba_0^{\xi}\vee\ba_1^{\xi}$ for all $\xi\in\Omega$,
there are a decomposition $\Omega=\bigcup\famm{\Omega_m}{m\in\gos}$ and
mappings $\bc_m\colon\Sem(m,\Omega_m)\times\Sem(m,\Omega_m)\to S$,
for $m\in\gos$, such that the following statements hold for every positive integer~$m$:
\begin{enumerate}
\item $p\leq q$ implies that $\bc_m(p,q)=0$, for all $p,q\in\Sem(m,\Omega_m)$;

\item $\bc_m(p,r)\leq\bc_m(p,q)\vee\bc_m(q,r)$, for all
$p,q,r\in\Sem(m,\Omega_m)$;

\item $\bc_m(p\vee q,r)=\bc_m(p,r)\vee\bc_m(q,r)$, for all
$p,q,r\in\Sem(m,\Omega_m)$;

\item $\bc_m(\one,\zero)=\be$;

\item The inequality $\bc_m((k+1)\cdot\dxi,k\cdot\dxi)\leq\ba_{\eps(k)}^{\xi}$ holds,
for all $\xi\in\Omega_m$ and all $k\leq m$.
\end{enumerate}
If, for a fixed $m\in\gos$, we can always take $\Omega_m=\Omega$ while $\Omega_n=\es$ for all $n\neq m$, we say that $S$ satisfies $\CLR_m(\be)$.
\end{definition}

The statement $\CLR(\be)$ is an analogue, for arbitrary lattices, of the
`uniform refinement property' introduced in \cite{WeURP}, denoted by
`$\mathrm{URP}^-$ at $\be$' in \cite{CLPSurv}.
It is easy to verify that for any \jzs s $S$ and $T$, any $\be\in S$, and any weakly distributive \jzh\ $\mu\colon S\to T$, if~$S$ satisfies $\CLR(\be)$, then~$T$ satisfies
$\CLR(\mu(\be))$. A similar observation applies to $\CLR_m$. Furthermore, a straightforward, although somewhat tedious, modification of the proof of Theorem~\ref{T:NonRepr}, gives, for example, the following result.

\begin{theorem}\label{T:CLR}
Let $L$ be a lattice and let $\be$ be a principal congruence of~$L$.
Then $\Conc L$ satisfies $\CLR(\be)$. Furthermore, if $L$ has $(m+1)$-permutable congruences \pup{where~$m$ is a given positive integer}, then $\Conc L$ satisfies $\CLR_m(\be)$. On the other hand, $\cG(\aleph_{\go+1})$ \pup{resp., $\cG(\aleph_{2^m})$} does not satisfy $\CLR(1)$ \pup{resp., $\CLR_m(1)$}.
\end{theorem}

\subsection{Open problems}\label{Su:pbs}
The most obvious problem suggested by the present paper is to fill the
cardinality gap between $\aleph_2$ and $\aleph_{\go}$. In the meantime, this problem has been solved by Pavel R\r{u}\v{z}i\v{c}ka \cite{Ruzi2}, who introduced a strengthening of Kuratowski's Free Set Theorem that made it possible to prove, by using the original Erosion Lemma (Lemma~\ref{L:erosion}) and modifications of both the Evaporation Lemma (Lemma~\ref{L:Evap}) and the Descent Lemma (Lemma~\ref{L:Descent}) the following result: \emph{For any set $\Omega$ such that $|\Omega|\geq\aleph_2$, there are no algebra~$L$ with a congruence-compatible structure of bounded semilattice and no weakly distributive \jzuh\ $\mu\colon\Conc L\to\cG(\Omega)$}. In fact, it is not hard to modify R\r{u}\v{z}i\v{c}ka's proof to establish that for $|\Omega|\geq\aleph_2$, the semilattice $\cG(\Omega)$ does not satisfy $\CLR(1)$ (cf. Subsection~\ref{Su:URP}).

The discussion in Subsection~\ref{Su:URP} about $\CLR$ and $\CLR_m$ also suggests the following problem.

\begin{problem}\label{Pb:CPE}
Prove that there exists a lattice $K$ such that for every positive integer~$m$, there is no lattice $L$ with $m$-permutable congruences such that\linebreak
$\Con K\cong\Con L$.
\end{problem}

Of course, it is sufficient to find a counterexample for each~$m$, as their direct product would then solve Problem~\ref{Pb:CPE}.

Now as we know that the answer to CLP is negative, a natural question
is the corresponding one for congruence-distributive varieties.

\begin{problem}\label{Pb:CDVar}
Is every algebraic distributive lattice isomorphic to the congruence
lattice of some algebra generating a congruence-distributive variety?
\end{problem}

Recall the classical open problem asking whether every algebraic
distributive lattice is isomorphic to the congruence lattice of some algebra
with finitely many operations. In view of Theorem~\ref{T:NonRepr}, we may
try to find the algebra with a \jzs\ (but not \jus) operation.

Kearnes proves in \cite{Kear05} that there exists an algebraic lattice
that is not isomorphic to the congruence lattice of any locally finite
algebra. In light of this result, the following question is natural.

\begin{problem}\label{Pb:LocFin}
Does there exist a lattice $L$ such that $\Con L$ is not isomorphic to
the congruence lattice of any locally finite lattice (resp., algebra)?
\end{problem}

In \cite{TuWe}, infinite semilattices considered
earlier in \cite{WeURP,PTW,NonMeas} are approximated by finite
semilattices, yielding, in particular, a $\set{0,1}^3$-indexed diagram of
finite Boolean semilattices that cannot be lifted, with respect to the
$\Conc$ functor, by congruence-permutable lattices. The methods used in
the present paper suggest that those works could be extended to find a
$\set{0,1}^{2^m+1}$-indexed diagram of finite Boolean semilattices that
cannot be lifted, with respect to the $\Conc$ functor, by lattices with $(m+1)$-permutable congruences.

T\r{u}ma and Wehrung prove in~\cite{Bowtie}
that there exists a diagram of finite Boolean semilattices, indexed by a
finite \emph{partially ordered set}, that cannot be lifted, with respect
to the $\Conc$ functor, by any diagram of lattices (or even algebras in
any variety satisfying a nontrivial congruence lattice identity). This
leaves open the following problem.

\begin{problem}\label{Pb:Anyfindiagr}
Prove that any diagram of finite distributive \jzs s and \jzh s, indexed
by a finite lattice, can be lifted, with respect to the $\Conc$ functor,
by a diagram of (finite?) lattices and lattice homomorphisms.
\end{problem}

We conclude with the following problem, which also appears, with a slightly different formulation, as \cite[Problem~10.6]{FCBook}.

\begin{problem}\label{Pb:ModLatt}
Prove that there exists a lattice $K$ such that there is no \emph{modular} lattice~$M$ with $\Con K\cong\Con M$.
\end{problem}

\section{Acknowledgment}
My deepest thanks go to the anonymous referees for their careful reading of the paper and valuable suggestions.

\section*{Added in proof}
A recent survey article partly devoted to CLP, written by George Gr\"atzer, just appeared, as ``Two Problems That Shaped a Century of Lattice Theory", Notices Amer. Math. Soc. \textbf{54}, no.~6 (2007), 696--707.

\end{document}